\documentclass{amsart}
\usepackage[utf8]{inputenc}
\usepackage[T1]{fontenc}
\usepackage[english,francais]{babel}
\usepackage{amsmath}
\usepackage{amssymb,amsthm}
\usepackage{amsfonts}

\newcommand{\Unc}{U^{\mathrm{nc}}}
\newcommand{\N}{\mathbb{N}}
\newcommand{\C}{\mathbb{C}}
\newcommand{\id}{\mathrm{id}}
\newcommand{\NC}{\mathrm{NC}}
\newcommand{\modu}[1]{\left\lvert #1 \right\rvert}
\newcommand{\fprod}{\text{\b{$\sqcup$}}}

\newtheorem{theorem}{Theorem}[section]
\newtheorem{lemma}[theorem]{Lemma}
\newtheorem{e-proposition}[theorem]{Proposition}

\newtheorem{e-definition}[theorem]{Definition\rm}

\newtheorem{theoreme}{Th\'eor\`eme}[section]

\newtheorem{definition}[theoreme]{D\'efinition\rm}

\title{Stochastic aspects of the unitary dual group}
\author{Isabelle Baraquin}
\email{isabelle.baraquin@univ-fcomte.fr}
\address{Laboratoire de math\'emathiques de Besan\c{c}on, UMR CNRS 6623, Universit\'e Bourgogne Franche-Comt\'e, 16 route de Gray, 25030 Besan\c{c}on cedex, France}

\begin{document}
\selectlanguage{english}
\begin{abstract}
\selectlanguage{english}

In this note we study asymptotic properties of the $*$-distribution of traces of some matrices, with respect to the free Haar trace on the unitary dual group. The considered matrices are powers of the unitary matrix generating the Brown algebra. We proceed in two steps, first computing the free cumulants of any R-cyclic family, then characterizing the asymptotic $*$-distributions of the traces of powers of the generating matrix, thanks to these free cumulants. In particular, we obtain that these traces are asymptotic $*$-free circular variables.
\vskip 0.5\baselineskip

\selectlanguage{francais}
\noindent{\bf R\'esum\'e}
\vskip 0.5\baselineskip
\noindent {\bf Aspects stochastiques du groupe dual unitaire}

Dans cette note, nous étudions la loi asymptotique de la trace de certaines matrices, par rapport \`a la trace de Haar libre sur le groupe dual unitaire. Ces matrices sont les puissances de la matrice unitaire qui engendre l'alg\`ebre de Brown. Nous proc\'edons en deux étapes. Tout d'abord, nous calculons les cumulants joints d'une famille de matrices R-cyclique. Nous caract\'erisons ensuite la $*$-distributions asymptotique des traces considérées, \`a l'aide des cumulants libres. En particulier, nous obtenons que ces traces sont des variables asymptotiquement circulaires et $*$-libres.
\paragraph*{Mots-cl\'es : } groupe dual unitaire, trace de Haar, probabilit\'es libres, cumulants, R-cyclicité, variables circulaires
\end{abstract}

\selectlanguage{english}
\maketitle

\medskip

\selectlanguage{francais}
\section*{Version fran\c{c}aise abr\'eg\'ee}
Dans cette note, nous étudions le comportement asymptotique des $\chi(u^p)$, traces des puissances de la matrice unitaire $u$ d\'efinissant l'alg\`ebre de Brown $\Unc_n$, par rapport à la trace de Haar libre $h$ d\'etermin\'ee par C\'ebron et Ulrich \cite{Michael}, et par McClanahan \cite{McClanahan}.

Pour cela, nous remarquons tout d'abord que la propriété de $h$ rappelée en Proposition \ref{cor2.8} entraîne que $\{u, u^*\}$ est une famille de matrices $R$-cyclique, c'est-à-dire que les cumulants libres $\kappa_r$ de coefficients $(u^e)_{ij}$ s'annulent dès lors que les indices ne sont pas cycliques : 
\begin{definition}[\cite{Shlya}]
Soit $\mathcal{A}$ une algèbre. Une famille de matrices de $\mathcal{M}_n(\mathcal{A})$, notée $\{A_l = (a_{ij}^{(l)})_{1 \leq i,j \leq n}\}_{1 \leq l \leq s}$, est appelée \emph{R-cyclique} si $\kappa_r(a_{i_1j_1}^{(l_1)}, \ldots, a_{i_rj_r}^{(l_r)}) = 0$ lorsqu'il n'est pas vrai que $j_1 = i_2$, $j_2 = i_3$, ..., $j_{r-1} = i_r$ et $j_r = i_1$.
\end{definition}

Nous pouvons donc nous ramener au calcul du cumlant libre des traces des puissances d'une famille de matrices R-cyclique, dont les cumulants considérés dans la définition ci-dessus ne dépendent pas des indices $i_1$, $\ldots$, $i_r$. L'application de ce résultat dans notre cas, nous permet de montrer que :
\begin{theoreme}
La famille de variables al\'eatoires $\left(\chi(u^p)\right)_{p \geq 1}$ est asymptotiquement une famille  $*$-libre de variables circulaires, d'esp\'erance $0$ et de covariance $1$.
\end{theoreme}

\selectlanguage{english}

\section{Introduction}

Diaconis, Shahshahani and Evans \cite{shah,evans} show that the traces of powers of a matrix chosen at random in the unitary (respectively orthogonal) group behave asymptotically like independent complex (resp. real) Gaussian random variables. Later, Banica, Curran and Speicher investigate the case of easy quantum groups in \cite{BCS2011}, and obtain similar results in the context of free probability for free orthogonal groups. 

In \cite{Michael}, C\'ebron and Ulrich study the Haar states according to the five notions of convolution (free, tensor, boolean, monotone and anti-monotone) of the unitary dual group $U\langle n \rangle$. They prove in particular that there is no Haar state for each of the five notions of convolution, and even no Haar trace for the boolean, monotone or anti-monotone convolution, for $n \geq 2$. They also define a faithful Haar trace on $U\langle n \rangle$ for the free convolution, denoted $h$, which is  in fact equal to the state given by McClanahan in \cite{McClanahan}.

The aim of this note is to extend the study of Diaconis, Shashahani and Evans to the framework of the unitary dual group and its free Haar trace. The paper is organized as follows. We first introduce the tools of our study. Then, in the last section, we discuss the computation of the joint cumulants of traces of powers of an R-cyclic family and determine the asymptotic $*$-distribution of the traces of powers of the generating matrix, with respect to the free Haar trace.

\section{Preliminaries}

We recall here some facts about the unitary dual group, free cumulants and R-cyclicity.

\subsection{The unitary dual group}
Let $n \geq 1$, and $\Unc_n$, sometimes called the Brown algebra, be the noncommutative $*$-algebra generated by $n^2$ elements $\{u_{ij}\}_{1 \leq i,j \leq n}$ such that the matrix $u = (u_{ij})_{1 \leq i,j\leq n}$ is unitary. It is possible to endow this algebra with a structure of dual group in the sense of Voiculescu \cite{voiculescu}, $U\langle n\rangle = (\Unc_n, \Delta, \delta, \Sigma)$, called the unitary dual group. This is a generalization of the notion of groups, like Hopf algebras, but using the free product instead of the tensor product.

\begin{e-definition}
Let $\mathcal{A}$ and $\mathcal{B}$ be unital $*$-algebras. The \emph{free product of $\mathcal{A}$ and $\mathcal{B}$} is the unique unital $*$-algebra $\mathcal{A} \sqcup \mathcal{B}$ with two $*$-homomorphisms $i_\mathcal{A} \colon \mathcal{A} \to \mathcal{A} \sqcup \mathcal{B}$ and $i_\mathcal{B} \colon \mathcal{B} \to \mathcal{A} \sqcup \mathcal{B}$, such that, for all $*$-homomorphisms $f \colon \mathcal{A} \to \mathcal{C}$ and $g \colon \mathcal{B} \to \mathcal{C}$, there exists a unique $*$-homomorphism $f \sqcup g \colon \mathcal{A} \sqcup \mathcal{B} \to \mathcal{C}$ satisfying $f = (f \sqcup g) \circ i_\mathcal{A}$ and $g = (f \sqcup g) \circ i_\mathcal{B}$.
\end{e-definition}

We sometimes refer to $\mathcal{A}$ and $\mathcal{B}$ as the left and right legs of the free product $\mathcal{A} \sqcup \mathcal{B}$. Therefore, for each $a \in \mathcal{A}$ and $b \in \mathcal{B}$, we denote $i_\mathcal{A}(a)$ and $i_\mathcal{B}(b)$ by $a^{(1)}$ and $b^{(2)}$, respectively. 

Let $f \colon \mathcal{A}_1 \to \mathcal{A}_2$ and $g \colon \mathcal{B}_1 \to \mathcal{B}_2$ be unital $*$-homomorphisms between the four unital $*$-algebras $\mathcal{A}_1$, $\mathcal{A}_2$, $\mathcal{B}_1$ and $\mathcal{B}_2$. Then we denote by $f \fprod g \colon \mathcal{A}_1 \sqcup \mathcal{B}_1 \to \mathcal{A}_2 \sqcup \mathcal{B}_2$ the unital $*$-homomorphism given by the free product $(i_{\mathcal{A}_2} \circ f) \sqcup (i_{\mathcal{B}_2} \circ g)$.

\begin{e-definition}
Let $n \geq 1$. The \emph{unitary dual group $U\langle n \rangle$} is defined by the unital $*$-algebra $U^{\mathrm{nc}}_n$ and three unital $*$-homomorphisms $\Delta \colon \Unc_n \to \Unc_n \sqcup \Unc_n$, $\delta \colon \Unc_n \to \C$ and $\Sigma \colon \Unc_n \to \Unc_n$, such that
\begin{itemize}
\item $\Unc_n$ is the Brown algebra, generated by the $u_{ij}$'s satisfying 
\[\forall 1\leq i,j \leq n, \ \sum_{k = 1}^n u_{ki}^*u_{kj} = \delta_{ij} = \sum_{k = 1}^n u_{ik}u_{jk}^*\text{ ,}\]
\item the map $\Delta$ is a coassociative coproduct, i.e.\ $(\id \fprod \Delta)\circ \Delta = (\Delta \fprod \id)\circ \Delta$, given on the generators by $\Delta(u_{ij}) = \sum\limits_{k = 1}^n u_{ik}^{(1)} u_{kj}^{(2)}$,
\item the map $\delta$ is a counit, i.e.\ $(\delta \fprod \id)\circ \Delta = \id = (\id \fprod \delta)\circ \Delta$, given by $\delta(u_{ij}) = \delta_{ij}$,
\item the map $\Sigma$ is a coinverse, i.e.\ $(\Sigma \sqcup \id)\circ \Delta = \delta(\cdot) 1_{\Unc_n} = (\id \sqcup \Sigma)\circ \Delta$, given by $\Sigma(u_{ij}) = u_{ji}^*$.
\end{itemize}
\end{e-definition}

\subsection{Free cumulants}

Note that the free cumulants characterize random variables and free independence. We will use them to study the asymptotic law of $\chi(u^p)$, as Banica, Curran and Speicher in \cite{BCS2011}. Let us introduce some notations:
\begin{itemize}
\item $[s]$ denotes the set $\{1, 2, \ldots, s\}$ for each $s \geq 1$,
\item the set of noncrossing partitions of $[s]$ is denoted $NC(s)$,
\item $NC_2(s)$ corresponds to the noncrossing pairings.
\end{itemize}

The family of free cumulants $(\kappa_r)_{r \in \N}$ is uniquely characterized by the corresponding multiplicative family of functionals satisfying the moment-cumulant formula \cite[equation (11.7)]{NicaSpeicher}:
\[\forall s \in \N, \forall \{a_i\}_{1 \leq i \leq s} \subset \Unc_n,\ h(a_1 \ldots a_s) = \sum_{\pi \in \NC(s)} \kappa_\pi[a_1, \ldots, a_s]\]
where $\kappa_\pi[a_1, \ldots, a_s] = \prod\limits_{V \in \pi} \kappa_V[a_1, \ldots, a_s] := \prod\limits_{\substack{V \in \pi\\V = \{v_1 < \ldots < v_l\}}} \kappa_l(a_{v_1}, \ldots, a_{v_l})$.

In particular, they satisfy the following property:
\begin{e-proposition}[{\cite[equation (11.11)]{NicaSpeicher}}]
\label{prod}
Let $s, r \in \N$ and $\underline{p} \in [n]^r$ be given such that the sum of the $p_i$'s is $s$. Then, for any $\{a_i\}_{1 \leq i \leq s} \subset \Unc_n$,
\[\kappa_r\left(a_1 \ldots a_{p_1}, a_{p_1+1} \ldots a_{p_1 + p_2}, \ldots, a_{s - p_r + 1} \ldots a_s\right) = \sum_{\substack{\pi \in NC(s)\\ \pi \vee \gamma_{\underline{p}} = 1_s}} \kappa_\pi[a_1, \ldots, a_s]\]
where $\gamma_{\underline{p}}$ denotes the noncrossing partition associated to the multi-index $\underline{p}$, i.e.\ $\gamma_{\underline{p}} = \left\{ \{1, \ldots, p_1\}, \ldots, \{s - p_r + 1, \ldots, s\}\right\}$.
\end{e-proposition}

C\'ebron and Ulrich \cite{Michael} compute the free cumulants associated to the free Haar trace $h$ of the generators of the Brown algebra $u_{ij}$ and their adjoints $(u^*)_{ij} = u^*_{ji}$.

\begin{e-proposition}[\cite{Michael}]
\label{cor2.8}
The free cumulants of $(u_{ij})_{1 \leq i,j \leq n}$ and $((u^*)_{ij})_{1 \leq i,j \leq n}$ in the noncommutative probability space $(\Unc_n, h)$ are given as follows.

Let $1 \leq i_1, j_1, \ldots, i_r, j_r \leq n$ and $\epsilon_1, \ldots, \epsilon_r$ be either $\emptyset$ or $*$. If the indices are cyclic (i.e. if $j_{l-1} = i_l$ for $2 \leq l \leq r$ and $i_1 = j_r$), $r$ is even and the $\epsilon_i$ are alternating we have
\[\kappa_r \left( (u^{\epsilon_1})_{i_1j_1}, \ldots, (u^{\epsilon_r})_{i_rj_r}\right) = n^{1-r} (-1)^{\frac{r}{2}-1} C_{\frac{r}{2}-1}\]
where $C_i = \frac{(2i)!}{(i+1)!i!}$ designate the Catalan numbers.

Otherwise, the left-hand side term is equal to zero.
\end{e-proposition}

\subsection{R-cyclicity}

\begin{e-definition}
For an algebra $\mathcal{A}$, a family of matrices $\{A_l = (a_{ij}^{(l)})_{1 \leq i,j \leq n}\}_{1 \leq l \leq s}$ in $\mathcal{M}_n(\mathcal{A})$ is called \emph{R-cyclic} if $\kappa_r(a_{i_1j_1}^{(l_1)}, \ldots, a_{i_rj_r}^{(l_r)}) = 0$ whenever it is not true that $j_1 = i_2$, $j_2 = i_3$, ..., $j_{r-1} = i_r$ and $j_r = i_1$.
\end{e-definition}

Thus Proposition \ref{cor2.8} ensures that $\{u, u^*\}$ is an R-cyclic family. Moreover, the cyclic free cumulants do not depend on the indices.

Since we want to look at powers of the generating matrix $u$, we need th following property:
\begin{e-proposition}[{\cite[Theorem 4.3]{Shlya}}]
\label{cyclic}
Let $(\mathcal{A}, \phi)$ be a noncommutative probability space. Let $d$ be a positive integer, and let $A_1, \ldots, A_s$ be an R-cyclic family of matrices in $\mathcal{M}_d(\mathcal{A})$. We denote by $\mathcal{D}$ the algebra of scalar diagonal matrices in $\mathcal{M}_d(\mathcal{A})$, and by $\mathcal{C}$ the subalgebra of $\mathcal{M}_d(\mathcal{A})$ which is generated by $\{A_1, \ldots, A_s\} \cup \mathcal{D}$. Then every finite family of matrices from $\mathcal{C}$ is R-cyclic.
\end{e-proposition}

In particular every finite subset of $\{u^k\}_{k \geq 1} \cup \{(u^*)^k\}_{k \geq 1}$ is also an R-cyclic family of matrices.

\section{Main result}

\begin{theorem}
The family of variables $\left(\chi(u^p)\right)_{p \geq 1}$ is asymptotically a family of $*$-free circular variables of mean $0$ and covariance $1$.
\end{theorem}

To prove this, let us use the R-cyclicity of $u$ and $u^*$. First, let us note that $\chi(u^p)^e = \chi((u^e)^p)$ for any $p \geq 1$ and $e \in \{\emptyset, *\}$. This means that we can see the calculation of $\kappa_s(\chi(u^{p_1})^{e_1}, \ldots, \chi(u^{p_s})^{e_s})$ in a more general framework and compute $\kappa_s(\chi(A_{l_1}^{p_1}), \ldots, \chi(A_{l_s}^{p_s}))$ for $\{A_l\}_{l \in I}$ an R-cyclic family of matrices.
\begin{lemma}
Let $\{A_l\}_{l \in I}$ be an R-cyclic family of matrices such that the cumulants $\kappa_\pi\left[a^{(l_1)}_{i_1 j_2}, \ldots, a^{(l_s)}_{i_s j_s}\right]$ depend only on the cyclicity of the indices. Let us denote by $\kappa_\pi\left[a^{(l_1)}, \ldots, a^{(l_s)}\right]$ the common value of the cumulants with cyclic indices, i.e.\ such that $j_1 = i_2$, $\ldots$, $j_s = i_1$. Then
\[\kappa_s(\chi(A_{l_1}^{p_1}), \ldots, \chi(A_{l_s}^{p_s})) = n^{p + 2 - s}\sum_{\substack{\pi \in NC(p)\\ \pi \vee \gamma_{\underline{p}} = 1_p}}  n^{-\modu{\pi}} \kappa_\pi\left[a^{(l_1)}, \ldots, a^{(l_1)}, a^{(l_2)}, \ldots, a^{(l_s)}\right] \text{ .}\]
\end{lemma}

Note that this calculation is similar to the one in the last section of \cite{Mingo}, and we will use similar arguments. Since $\{A_{l_1}^{p_1}, \ldots, A_{l_s}^{p_s}\}$ is also an R-cyclic family, by Proposition \ref{cyclic}, if $p$ is the sum of all the $p_i$'s and $a_{ij}^{(l)}$ denotes the coefficient $(i,j)$ of the matrix $A_l$,
\[\kappa_s(\chi(A_{l_1}^{p_1}), \ldots, \chi(A_{l_s}^{p_s})) = \hspace{-30pt}\sum_{\substack{1 \leq i_1, \ldots, i_p \leq n\\ i_1 = i_{p_1 + 1} = \ldots = i_{p - p_s +1}}}\hspace{-30pt}\kappa_s\left(a^{(l_1)}_{i_1 i_2} \ldots a^{(l_1)}_{i_{p_1} i_{p_1+1}}, \ldots, a^{(l_s)}_{i_{p - p_s + 1} i_{p - p_s + 2}} \ldots a^{(l_s)}_{i_p i_1}\right) \text{ .}\]
By Proposition \ref{prod}, this is equal to
\[\sum_{\substack{1 \leq i_1, \ldots, i_p \leq n\\ i_1 = i_{p_1 + 1} = \ldots = i_{p - p_s +1}}} \sum_{\substack{\pi \in NC(p)\\ \pi \vee \gamma_{\underline{p}} = 1_p}} \kappa_\pi\left[a^{(l_1)}_{i_1 i_2}, \ldots, a^{(l_1)}_{i_{p_1} i_{p_1+1}}, \ldots, a^{(l_s)}_{i_p i_1}\right] \text{ .}\]
By definition of $\kappa_\pi$, we can restrict ourselves to the study of a block $V$ of $\pi$ and look at $\kappa_r(a^{(\lambda_{v_1})}_{i_{v_1}i_{v_1+1}}, \ldots, a^{(\lambda_{v_r})}_{i_{v_r}i_{v_r+1}})$ with $\underline{\lambda} = (l_1, \ldots, l_1, l_2, \ldots, l_s) \in I^p$ and where $V$ is the block $\{v_1 < \ldots < v_r\}$. In order to have a non zero contribution, the indices have to be cyclic, i.e.\ to satisfy
\[\forall 1 \leq j \leq r-1, i_{v_j+1} = i_{v_{j+1}} \text{ and } i_{v_r + 1} = i_{v_1} \text{ .}\]
Let us denote by $\sigma_\pi$ the permutation associated to the partition $\pi$, by considering the elements of a block of $\pi$ in increasing order as a cycle of $\sigma_\pi$. Hence the conditions above can be written as $i_{\gamma(v_i)} = i_{\sigma_\pi(v_i)}$ where $\gamma = (1, 2, \ldots, p)$. Since this should be true for each block of $\pi$, it means that $i_j = i_{\gamma \circ \sigma_\pi^{-1}(j)}$ for all $1 \leq j \leq p$. Thus, we get
\[\kappa_s(\chi(A_{l_1}^{p_1}), \ldots, \chi(A_{l_s}^{p_s})) =\hspace{-5pt}\sum_{\substack{\pi \in NC(p)\\ \pi \vee \gamma_{\underline{p}} = 1_p}} \hspace{2pt} \sum_{\substack{1 \leq i_1, \ldots, i_p \leq n\\ i_1 = i_{p_1 + 1} = \ldots = i_{p - p_s +1}\\i_j = i_{\gamma \circ \sigma_\pi^{-1}(j)}}} \hspace{-20pt} \kappa_\pi\left[a^{(l_1)}_{i_1 i_2}, \ldots, a^{(l_1)}_{i_{p_1} i_{p_1+1}}, \ldots, a^{(l_s)}_{i_p i_1}\right] \text{ .}\]

Since, moreover $\kappa_\pi\left[a^{(l_1)}_{i_1 i_2}, \ldots, a^{(l_s)}_{i_p i_1}\right]$ depends only on the cyclicity of the indices, and is denoted by $\kappa_\pi\left[a^{(l_1)}, \ldots, a^{(l_1)}, \ldots, a^{(l_s)}\right]$, we obtain
\[\kappa_s(\chi(A_{l_1}^{p_1}), \ldots, \chi(A_{l_s}^{p_s})) = \sum_{\substack{\pi \in NC(p)\\ \pi \vee \gamma_{\underline{p}} = 1_p}}  \kappa_\pi\left[a^{(l_1)}, \ldots, a^{(l_1)}, a^{(l_2)}, \ldots, a^{(l_s)}\right]c_{\pi}\]
where $c_\pi$ denotes the quantity
\[\#\left\{\underline{i} \in \{1, \ldots, n\}^p, \; i_1 = i_{p_1 + 1} = \ldots = i_{p - p_s +1}, \ \forall 1\leq j \leq p, i_j = i_{\gamma \circ \sigma_\pi^{-1}(j)} \right\} \text{ .}\]

Thanks to \cite[Lemma 14]{Tan}, $\pi \vee \gamma_{\underline{p}} = 1_p$ if and only if $\sigma_{\pi}^{-1} \circ \gamma$ separates $p_1$, $p_1 + p_2$, $\ldots$ and $p$, which is equivalent to the fact that $1$, $p_1 + 1$, $\ldots$ and $p - p_s + 1$ are all in different blocks of $\gamma \circ \sigma_\pi^{-1}$. Thus, $c_\pi = n^{\#(\gamma \circ \sigma_\pi^{-1}) - s + 1}$, where $\#\sigma$ denotes the number of cycles in the cycle decomposition of a permutation $\sigma \in \mathfrak{S}_p$. Notes that $\gamma \circ \sigma_\pi^{-1}$ corresponds to the Kreweras complement \cite{NicaSpeicher} of $\pi$, denoted $K(\pi)$, conjugated by $\gamma$. Hence
\[\#(\gamma \circ \sigma_\pi^{-1}) = \modu{K(\pi)} = p + 1 - \modu{\pi}\]
where $\modu{\pi}$ is the number of blocks of $\pi$, and then $c_\pi = n^{p + 2 -s - \modu{\pi}}$, which proves the lemma.

In particular, in the dual unitary group endowed with the free Haar trace, we have $I = \{\emptyset, *\}$ and Proposition \ref{cor2.8} ensures that, if the partition $\pi$ is $\underline{\lambda}$-adapted with $\underline{\lambda} = (l_1, \ldots, l_1, l_2, \ldots, l_s)$,
\[\kappa_\pi\left[a^{(\lambda_1)}, \ldots, a^{(\lambda_p)}\right] = n^{\modu{\pi} - p} (-1)^{\frac{p}{2} - \modu{\pi}} \prod_{V \in \pi} C_{\frac{\#V}{2}-1}\]
otherwise the cumulant vanishes. Here, $\pi \in \NC(p)$ is said to be $\underline{\lambda}$-adapted when the following conditions are true for each block $V = \{v_1 < \ldots < v_l\}$ of the noncrossing partition $\pi$:
\begin{itemize}
\item $\#V := l \in 2\N$,
\item $\forall 1 \leq i \leq l-1, \ \lambda_{v_i} \neq \lambda_{v_{i+1}}$.
\end{itemize}

Finally, we get, with $\underline{\epsilon} = (e_1, \ldots, e_1, e_2, \ldots, e_s)$,
\[\kappa_s\left(\chi(u^{p_1})^{e_1}, \ldots, \chi(u^{p_s})^{e_s}\right) = n^{2-s} (-1)^{\frac{p}{2}} \sum_{\substack{\pi \in NC(p)\\ \pi \vee \gamma_{\underline{p}} = 1_p\\ \underline{\epsilon}-\text{adapted}}}\hspace{-10pt}(-1)^{\modu{\pi}} \prod_{V \in \pi} C_{\frac{\#V}{2}-1} \text{ .}\]

If $s > 2$, it is clear that the cumulants vanish asymptotically. If $s = 2$, the cumulant is non zero if and only if $p_1 = p_2$ and $e_1 \neq e_2$, in this case the cumulant equals $1$. Moreover, if $s = 1$, there is no $(e)$-adapted partition, and the cumulant is zero. This is the description of the cumulants of a $*$-free family of circular random variables.

\section*{Acknowledgements}
This work was supported by the French "Investissements d'Avenir" program, project ISITE-BFC (contract ANR-15-IDEX-03).


\begin{thebibliography}{03}
\bibitem{BCS2011}
T. Banica, S. Curran and R. Speicher, Stochastic aspects of easy quantum groups, {Probab. Theory Related Fields} 149 (2011) 435--462.

\bibitem{Michael}
G. C\'ebron and M. Ulrich, Haar states and L\'evy processes on the unitary dual group, {J. Funct. Anal.} 270 (2016) 2769--2811.

\bibitem{evans}
P. Diaconis and S. N. Evans, Linear functionals of eigenvalues of random matrices, {Trans. Amer. Math. Soc.} 353 (2001) 2615--2633.

\bibitem{shah}
P. Diaconis and M. Shahshahani, On the eigenvalues of random matrices, {J. Appl. Probab.} 31 (1994) 49--62.

\bibitem{McClanahan}
K. McClanahan, {{$C^*$}-algebras generated by elements of a unitary matrix}, {J. Funct. Anal.} 107 (1992) 439--457

\bibitem{Mingo}
J.~A. Mingo and R. Speicher, Schwinger-{D}yson equations: classical and quantum, {Probab. Math. Statist.}, 33 (2013) 275--285.

\bibitem{Tan}
J.~A. Mingo, R. Speicher, and E. Tan, Second order cumulants of products, {Trans. Amer. Math. Soc.}, 361 (2009) 4751--4781.

\bibitem{NicaSpeicher}
A. Nica and R. Speicher, {Lectures on the combinatorics of free probability}, Lecture Note Ser., vol. 335, Cambridge University Press, 2006.

\bibitem{Shlya}
A. Nica, D. Shlyakhtenko, and R. Speicher, {$R$}-cyclic families of matrices in free probability, {J. Funct. Anal.}, 188 (2002) 227--271.

\bibitem{voiculescu}
D. Voiculescu, Dual algebraic structures on operator algebras related to free products, {J. Operator Theory} 17 (1987) 85--98.
\end{thebibliography}
\end{document}